\documentclass[smallextended,referee,envcountsect]{svjour3}
\smartqed
\usepackage{graphicx}
\usepackage{amssymb}
\usepackage{color}
\begin{document}

\title{Characterization of the norm-based robust solutions in vector optimization}

\subtitle{}

\author{Morteza Rahimi \and Majid Soleimani-damaneh}
\authorrunning{M. Rahimi \and M. Soleimani-damaneh}

\institute{\baselineskip=0.5\baselineskip
Morteza Rahimi\at
School of Mathematics, Statistics and Computer Science, College of Science, University of
               Tehran, Tehran, Iran; E-mail: rahimi.math@ut.ac.ir
\and
Majid Soleimani-damaneh (Corresponding author)\at
School of Mathematics, Statistics and Computer Science, College of Science, University of
               Tehran, Tehran, Iran; E-mail: soleimani@khayam.ut.ac.ir
}

\date{\\ Received: date / Accepted: date}

\maketitle

\begin{abstract}
In this paper, we study the norm-based robust (efficient) solutions of a Vector Optimization Problem (VOP). We define two kinds of non-ascent directions in terms of Clarke's generalized gradient and characterize norm-based robustness by means of the newly-defined directions. This is done under a basic Constraint Qualification (CQ). We extend the provided characterization to VOPs with conic constraints. Moreover, we derive a necessary condition for norm-based robustness utilizing a nonsmooth gap function.
\end{abstract}


\section{Introduction}

In optimization models arisen in practice, the Decision Maker (DM)/manager/ user is often faced with uncertainty. Robust optimization, as one of the leading  tools for dealing with uncertainty, has been the subject of many publications in recent decades; see, e.g., \cite{ben-1,ben-2,ber,ehr-2,fli-2,geo,gob-1,ide-1,ide-3,kur,zam} among others.

In the current work, we concentrate on robustness in Vector Optimization Problems (VOPs). We are going to investigate the efficient solutions which are insensitive against small perturbation in objective function data. In the following, we briefly review some relevant works. To this end, we classify the existing robustness notions to three classes: worst-case, set-based, and norm-based.

Worst-case robustness in multi-objective programming has been studied by Ehrgott et al. \cite{ehr-2}, Fliege and Werner \cite{fli-2}, and Kuroiwa and Lee \cite{kur}. This kind of robustness deals with a conservative over-estimator of the function on the whole uncertainty set \cite{fli-2}.

Set-based robustness has been appearing in some recent works by Ehrgott et al. \cite{ehr-2}, Ide and K$\ddot{o}$bis \cite{ide-1}, and Ide and Sch$\ddot{o}$bel \cite{ide-3}. The main idea behind set-based robustness is to compare the objective function values, taking the whole uncertainty set into account, by means of the set relations.

Norm-based robustness has been introduced by Georgiev et al. \cite{geo}  and then has been developed by Goberna et al. \cite{gob-1} and Zamani et al. \cite{zam}. This notion, which is useful for modelling in an unbalanced situation, refers to the efficient solutions which remain efficient under small perturbations.

In a very recent work, Rahimi and Soleimani-damaneh \cite{our-new} have defined and investigated robust efficient solutions of a nonlinear VOP. To the best of our knowledge, this is the only work in the literature dealing with robustness in vector optimization. In \cite{our-new}, we have defined, compared and characterized robust solutions from various standpoints. Furthermore, we have studied the connections between norm-based robust efficiency, strict efficiency, isolated efficiency, and proper efficiency.

In the current work, we are going to define two kinds of non-ascent directions in terms of Clarke's generalized gradient and then characterize norm-based robustness with respect to the newly-defined directions. To this end, we apply an appropriate Constraint Qualification (CQ). We derive various necessary and sufficient conditions for norm-based robustness.

The rest of the paper is organized as follows. In Section \ref{prelim}, we provide the required preliminaries. In Section \ref{full}, a characterization of the norm-based robust efficient solutions, in terms of tangent/normal cone and aforementioned directions, is given. Section \ref{conic} is devoted to investigation of the problem for VOPs with conic constraints. Section \ref{gap} concludes the paper by studying robustness invoking a new nonsmooth gap function.

\section{Preliminaries}\label{prelim}

This section contains some preliminaries which are used in the rest of the work. Given $x,y \in \mathbb{R}^n$, two notations $x^T$ and $x^Ty$ stand for the transpose of $x$ and the inner product of $x,y$, respectively. We denote the convex hull, the interior, and the closure of a nonempty set $\Omega \subseteq \mathbb{R}^n$ by $co\,\Omega$, $int\,\Omega$,
and $cl\,\Omega$, respectively.

Given $K \subseteq \mathbb{R}^p$
is said to be a cone if $x \in K$ and $\lambda \geq 0$ imply $\lambda x \in K$.
A cone $K$ is called pointed if $K \cap -K =\lbrace 0 \rbrace$; and it is said to be an ordering cone if it is nontrivial, closed, pointed and convex.
For instance, $$\mathbb{R}^p_+ := \{ d\in \mathbb{R}^p~:~d_i\geq 0,~ i=1,2,\ldots ,p\},$$ is called the natural ordering cone.

Given an infinite index set $J$, we set
$$\mathbb{R}^{|J|}_\geq :=\bigg\{ \mu :J\longrightarrow \mathbb{R}_\geq : card\{ j\in J \, |\,  \mu_j :=\mu(j)>0\} <\infty \bigg\}.$$

The cone and the convex cone generated by $\Omega \subseteq \mathbb{R}^n$ are denoted by $cone(\Omega)$ and
$pos(\Omega)$, respectively.
Indeed, $$cone(\Omega):=\{ \alpha y:~ y \in \Omega ,~ \alpha \geq 0\},$$
$$pos(\Omega):=\bigg\{\sum_{i=1}^m \lambda_iy_i ~ :~ m\in \mathbb{N},~\lambda_i\geq 0,~y_i\in
\Omega,~i=1,2,\ldots,m\bigg\}.$$

Given an infinite index set $T$ and a collection of nonempty convex sets $\{\Omega_t \subseteq \mathbb{R}^n :t\in T\}$, we have $pos\bigg(\displaystyle\bigcup_{t\in T} \Omega_t\bigg)=\bigcup_{\hat{T} \in \Sigma} pos\bigg(\displaystyle\bigcup_{t\in \hat{T}} \Omega_t\bigg),$
where $\Sigma$ is the set of all nonempty finite subsets of $T$; see \cite[Theorem 3.3]{roc}.



The polar cone and the strict polar cone of a set $K \subseteq \mathbb{R}^n$, denoted by $K^*$ and $K^{*\circ}$, respectively, are defined as
$$K^*:=\big\{y\in \mathbb{R}^p:~ y^Tx\leq 0,~\forall x\in K \big\},$$
$$K^{*\circ}:=\big\{y\in \mathbb{R}^p:~ y^Tx< 0,~\forall x\in K\setminus \{0\} \big\}.$$
The tangent cone to $\Omega$ at $\bar{x} \in \Omega$, denoted by
$T_{\Omega}(\bar x)$, is defined as
$$T_{\Omega}(\bar x):=\biggl\{d\in \mathbb{R}^n~:~\exists\big(\{x_\nu\}_\nu\subseteq \Omega,~\{t_\nu\}_\nu\subseteq \mathbb{R}\big);~t_\nu\downarrow 0,
~\frac{x_\nu-\bar x}{t_\nu}\longrightarrow d\biggl\}.$$
The normal cone to $\Omega$ at $\bar{x} \in \Omega$, denoted by $N_{\Omega}(\bar x)$, is defined as polar of the tangent cone, i.e., $$N_{\Omega}(\bar x):={T_{\Omega}(\bar x)}^*.$$

We use the Euclidean norm, i.e., $\Vert d\Vert =\sqrt{d^Td}$, and set $$\mathbb{S} := \{x \in \mathbb{R}^n: \Vert x\Vert = 1\}.$$

Let $K\subseteq \mathbb{R}^p$ be an ordering cone. A vector-valued function $f:\mathbb{R}^n \rightarrow\mathbb{R}^p$
is called $K$-convex if for any $x,y\in \mathbb{R}^n$ and $\lambda \in [0,1]$,
$$f(\lambda x+ (1-\lambda)y) -\lambda f(x)-(1-\lambda) f(y) \in -K.$$

The classic and Clarke's generalized directional derivatives \cite{cla} are used in the presence of nonsmooth data.

\begin{definition}
Let $h: \mathbb{R}^n\rightarrow \mathbb{R}$ and $\bar x,d\in \mathbb{R}^n$ be given. The directional derivative of $h$ at $\bar x$ in the direction $d$, denoted by $h'(\bar x; d)$, is defined as $$h'(\bar x;d):=
\displaystyle\lim_{t\downarrow 0}\frac{h(\bar x+td)-h(\bar x)}{t}.$$
\end{definition}
\begin{definition}
Let $h: \mathbb{R}^n\rightarrow \mathbb{R}$ be convex. The set of all subgradients of $h$
at $\bar x$, denoted by $\partial h(\bar x)$, is defined as
$$\partial h(\bar x)
:= \{ \zeta\in \mathbb{R}^n : h(x)-h(\bar x)\geq
\zeta^T(x-\bar x),~~\forall x\in \mathbb{R}^n\}.$$\end{definition}

\begin{definition}
\cite{cla} Let $h: \mathbb{R}^n\rightarrow \mathbb{R}$ be locally Lipschitz at $\bar x\in \mathbb{R}^n$. The Clarke's generalized directional
derivative of $h$ at $\bar x$ in the direction $d\in\mathbb{R}^n$, denoted by $h^\circ
(\bar x; d)$, is defined as $$h^\circ (\bar x; d):=
\displaystyle\limsup_{\stackrel{x\rightarrow \bar x}{t\downarrow
0}}\frac{h(x+td) - h(x)}{t}.$$
\end{definition}

\begin{definition}
\cite[Definition 10.3]{cla}  Let $h: \mathbb{R}^n\rightarrow \mathbb{R}$ be locally Lipschitz at
$\bar x\in \mathbb{R}^n$. The Clarke's generalized gradient of $h$
at $\bar x$, denoted by $\partial_{c} h(\bar x)$, is defined as
$$\partial_{c} h(\bar x)
:= \{ \zeta\in \mathbb{R}^n : h^\circ(\bar x; d)\geq \zeta^Td,~~
\forall d\in \mathbb{R}^n\}.$$
\end{definition}

If $h: \mathbb{R}^n\rightarrow \mathbb{R}$ is convex, then $h^\circ (\bar x;\cdot)=h'(\bar x;\cdot)$ and $\partial_{c} h(\cdot)=\partial h(\cdot)$; see~\cite{cla}.

A function $h: \mathbb{R}^n \rightarrow \mathbb{R}$ is called regular at $\bar{x} \in \mathbb{R}^n,$ if it is locally Lipschitz at $\bar{x}$, and $h'(\bar{x}; d)$ exists satisfying $h^\circ(\bar x;d)=h'(\bar{x};d)$ for any $d\in \mathbb{R}^n$ \cite[Definition 10.12]{cla}. The set of regular functions contains that of convex functions \cite{cla}.

Consider the VOP,
\begin{equation}\label{prb}
\min\, f(x)~~s.t.~~x \in \Omega,
\end{equation}
in which $f:\Omega\subseteq\mathbb{R}^n\rightarrow \mathbb{R}^p$ is a vector-valued function with $p\geq 2$. Indeed, $$f(x)=(f_1(x), f_2(x),
\ldots, f_p(x))^T,~~x\in \Omega.$$ Here, $\Omega$ and $f$ are the feasible set and the objective function, respectively. Throughout the paper, we suppose $f_i$, $i=1,\ldots ,p$, are locally Lipschitz. 
 Also, we consider an ordering cone $K\subseteq \mathbb{R}^p$ with nonempty interior.


\begin{definition}\label{def1}
The vector $\bar{x}\in \Omega$ is called an efficient solution of
(\ref{prb}) w.r.t. $K$ if there exists no $x\in \Omega$ such that
$f(x)- f(\bar{x})\in -K\setminus \{0\}$.
\end{definition}

We close this section by definition of norm-based robust efficient solution for VOPs. This notion, introduced by Rahimi and Soleimani-damaneh \cite{our-new}, generalizes the concepts scrutinized by Georgiev et al. \cite{geo} and Zamani et al. \cite{zam}. Before going to the definition, we need some notations.  For an $m \times n$ matrix $C = [c_{ij}]$, the Frobenius
norm is defined as $$\Vert C\Vert =\left( \sum_{i,j} \vert c_{ij}\vert^2\right)^{1/2}.$$ The set of all real $m \times n$ matrices
is denoted by $M(m, n)$; and the
set of all matrices $C \in M(m, n)$ with $\Vert C\Vert < r$ is denoted by $M(m, n; r)$. Given $W \subseteq M(m, n)$ and $V \subseteq \mathbb{R}^n$, we define $WV := \{wv ~:~ w \in W,~ v \in V \}$.

\begin{definition}\label{rob}\cite{our-new}
 The vector $\bar{x}\in \Omega$ is called a norm-based robust efficient solution of (\ref{prb}) w.r.t. $K$, if there exists some scalar $r>0$ such that for any $C \in M(p,n;r)$, the vector $\bar{x}$ is an efficient solution of
\begin{equation}\label{prb2}
\min \ f(x)+Cx~~s.t.~~x\in\Omega,
\end{equation}
w.r.t. $K$. The scalar $r$ is called a robustness radius for $\bar x.$
\end{definition}

\section{Characterization}\label{full}

In this section, we provide a full characterization of norm-based robust efficient solutions, for VOPs, in terms of Clarke's generalized gradient. To this end, we define two kinds of non-ascent directions of the objective
functions. Notice that as $f:\mathbb{R}^n\longrightarrow \mathbb{R}^p$ is a vector-valued function, the members of $\partial_c f(\bar{x})$ are $n\times p$ matrices whose columns are $\partial_c f_i(\bar{x}),~ i=1,2, \ldots ,p$.

\begin{definition}
A vector $d\in \mathbb{R}^n$ is called a first kind non-ascent direction of $f$ at
$\bar{x}$ if $d^T\xi \mu^*\leq 0,$ for each $\xi\in\partial_c
f(\bar{x})$ and each $\mu^* \in -K^*$.
\end{definition}

\begin{definition}
A vector  $d\in \mathbb{R}^n$ is called a second kind non-ascent direction of $f$ at
$\bar{x}$ if $d^T\eta \leq 0,$ for each $\eta \in\partial_c
(\mu^*\circ f)(\bar{x})$ and each $\mu^* \in -K^*$.
\end{definition}

Hereafter, $G_1(\bar{x})$ and $G_2(\bar{x})$  denote the set of all first and second kind non-ascent directions of $f$
at $\bar x$, respectively. Due to the properties of polar cone, $d\in G_1(\bar x)$ implies $\xi^Td\in(-K)^{**}=-K$ for each $\xi\in\partial_c
f(\bar{x})$. Here, $(-K)^{**}$ stands for the polar of $(-K)^*.$ According to \cite[Proposition 10.15]{cla}, $G_1(\bar{x}) \subseteq G_2(\bar{x})$ is always true, and it holds as equality if $K=\mathbb{R}^p_+$ and $f_i$, $i=1,2,\ldots ,p$, are regular at $\bar{x}$ in the sense of Clarke.

\begin{definition}\label{cq1}
We say that Constraint Qualification 1 (CQ1) holds at
$\bar{x}$ if
$G_1(\bar{x})=G_2(\bar{x}).$
\end{definition}

It can be seen that CQ1 holds $\bar{x}$ if for any $\mu^* \in -K^*$,
\begin{equation}\label{condition}
\partial_c(\mu^* \circ f)(\bar{x})=\left\{A\mu^*: A\in\partial_cf(\bar{x})\right\}.
\end{equation}


Theorem \ref{3part} is one of the most important results of the paper.

\begin{theorem}\label{3part} Let $\bar{x} \in \Omega$.
\begin{itemize}
\item[(i)] If $\bar{x}$ is a norm-based robust efficient solution of
(\ref{prb}) w.r.t. $K$, then $$T_{\Omega}(\bar{x})\cap G_1(\bar{x})=\{0\}.$$
\item[(ii)] Let $\Omega$ be convex and $f$
be $K$-convex. If $T_{\Omega}(\bar{x})\cap G_2(\bar{x})=\{0\},$ then $\bar{x}$ is a norm-based robust efficient solution of
(\ref{prb}) w.r.t. $K$.
\item[(iii)] Let $\Omega$ be convex, $f$
be $K$-convex, and CQ1 hold at $\bar x$. Then, $\bar{x}$ is a norm-based robust efficient solution of
(\ref{prb}) w.r.t. $K$ if and only if $$T_{\Omega}(\bar{x})\cap G_1(\bar{x})=\{0\}.$$
\end{itemize}
\end{theorem}
{\it \textbf{Proof}}
(i) By indirect proof, assume that there exists a nonzero vector $d$ such that $d\in T_{\Omega}(\bar{x}) \cap G_1(\bar{x})$. Then,
there are two sequences $\{x_\nu\}\subseteq \Omega$ and $t_\nu \downarrow 0$ such that
$\frac{x_\nu-\bar{x}}{t_\nu}\longrightarrow d$ as $\nu\longrightarrow \infty$.
By \cite[Theorem 10.17]{cla}, for each $\nu$,
\begin{equation}\label{in2}
f(x_\nu)=f(\bar{x})+\xi_\nu^T(x_\nu-\bar{x}),
\end{equation}
where $\xi_\nu \in M(n,p)$; and the $i$th column of $\xi_\nu$ belongs
to $\partial_c f_i(y_\nu^i)$ for some $y_\nu^i\in(\bar{x},x_\nu)$.
Locally Lipschitzness of $f$ at $\bar{x}$ and $y_\nu^i\longrightarrow \bar x$ imply that the sequence $\{\xi_\nu\}$ is bounded and without loss of generality is convergent to some $\xi\in \partial_c f(\bar{x})$. So, due to $d\in G_1(\bar{x})$, we have
$\xi^T d \in -K.$ Now, assume that $r>0$ is a robustness radius of $\bar{x}$. As $intK \neq \emptyset$ and $d\neq 0$, there exists some $\bar C \in M(p,n;r)$ with $\bar Cd \in -int K.$ So, $\xi^Td +\bar{C}d \in -intK$ and
for sufficiently large $\nu$,
$$\xi_\nu^T (\frac{x_\nu-\bar{x}}{t_\nu}) +\bar{C} (\frac{x_\nu-\bar{x}}{t_\nu}) \in -intK \subseteq -K \setminus \{0\}.$$
This leads to
$$\xi_\nu^T (x_\nu-\bar{x}) +\bar{C} (x_\nu-\bar{x}) \in -K\setminus \{0\}.$$
Hence, according to (\ref{in2}),
$$f(x_\nu)+\bar{C} x_\nu -f(\bar{x}) -\bar C\bar{x} \in -K\setminus \{0\}.$$
This contradiction completes the proof of part (i).

(ii) By indirect proof, assume that there exist two sequences $\{C_\nu\} \subseteq M(p,n)$ and $\{x_\nu\}\subseteq\Omega$ such that $C_\nu\longrightarrow 0$ and for any $\nu\in \mathbb{N},$
\begin{equation}\label{a89n}
f(x_\nu)+C_\nu x_\nu- f(\bar{x})-C_\nu\bar{x} \in -K\setminus \{0\}.
\end{equation}
This implies
\begin{equation}\label{ma89}
(\mu^*\circ f)(x_\nu) -(\mu^*\circ f)(\bar{x}) + \mu^{*^T} C_\nu(x_\nu -\bar{x}) \leq 0, ~~~\forall \Big(\mu^* \in -K^*,~\nu\in\mathbb{N}\Big).
\end{equation}
Without loss of generality, assume that $d_\nu:=\frac{x_\nu-\bar{x}}{\|x_\nu-\bar{x}\|}$ converges to some $d\in \mathbb{R}^n$ with $\|d\| =1.$

Two cases for the sequence $\{x_\nu\}$ may occur; either it has a
subsequence convergent to $\bar{x}$ or it does not have any
subsequence convergent to $\bar{x}$. In the first case, without loss of generality, assume
$x_\nu\longrightarrow \bar{x}$.
Then, $$d=\displaystyle\lim_{\nu\to\infty}\frac{x_\nu-\bar{x}}{\|x_\nu-\bar{x}\|}\in T_\Omega(\bar{x}).$$ On the other hand, as $f$ is $K$-convex, for any $\nu\in\mathbb{N}$,
\begin{equation}\label{ma90}
\eta^T (x_\nu -\bar{x}) \leq (\mu^*\circ f)(x_\nu) -(\mu^*\circ f)(\bar{x}),~~~\forall \Big(\mu^* \in -K^* ,~ \eta \in \partial_c(\mu^*\circ f)(\bar{x})\Big).
\end{equation}
Combining (\ref{ma89}) and (\ref{ma90}) leads to
$$\begin{array}{ll}
&\eta^T (x_\nu -\bar{x}) +\mu^{*^T} C_\nu(x_\nu -\bar{x}) \leq 0,~~~~~~\forall \nu\in\mathbb{N}\vspace*{2mm}\\
\Longrightarrow &\eta^T \frac{x_\nu -\bar{x}}{\|x_\nu -\bar{x}\|} +\mu^{*^T} C_\nu \frac{x_\nu -\bar{x}}{\|x_\nu -\bar{x}\|} \leq 0,~~~\forall \nu\in\mathbb{N}\vspace*{2mm}\\
\stackrel{\nu\to \infty}{\Longrightarrow} &\eta^T d \leq 0,
\end{array}$$
 for any $\mu^* \in -K^*$ and any $\eta \in \partial(\mu^*\circ f)(\bar{x})$. This implies that $d\in G_2(\bar{x})$.
Therefore, $0\neq d\in T_\Omega (\bar{x}) \cap G_2(\bar{x}).$ This contradicts the assumption.

Now, we consider the second case: $\{x_\nu\}$ does not have any
subsequence convergent to $\bar{x}$. Therefore, without loss of
generality, there exists some scalar $r>0$ such that $\|x_\nu-\bar x\|>r$ for any $\nu \in \mathbb{N}$. By setting $x_\nu^{'}=\bar x+\frac{r}{\nu}d_\nu$, as $\Omega$ is convex, we have
$$x_\nu^{'}=\frac{r}{\nu\|x_\nu-\bar{x} \|} x_\nu +(1-\frac{r}{\nu\|x_\nu-\bar{x} \|} ) \bar{x}\in \Omega$$ for any $\nu\in \mathbb{N}.$ Furthermore, $x_\nu^{'}\to\bar x$. Hence, $$d=\displaystyle\lim_{\nu\to\infty}\frac{x_\nu^{'}-\bar x}{r/\nu}\in T_\Omega (\bar{x}).$$ Let $t\in(0,r)$ be arbitrary. Similar to above, $\bar x+td_\nu\in\Omega,$ for each $\nu$. On the other hand, from the $K$-convexity and the locally Lipschitzness of $f$,
$$\begin{array}{ll}
&f(\bar x+td_\nu)-\frac{t}{\|x_\nu-\bar{x} \|} f(x_\nu) - (1-\frac{t}{\|x_\nu-\bar{x} \|}) f(\bar{x}) \in -K\vspace{2mm} \\
\hspace*{1.5mm}\Longrightarrow & \frac{f(\bar x+td_\nu)-f(\bar{x})}{t}-\frac{1}{\|x_\nu-\bar{x}\|} (f(x_\nu)-f(\bar{x})) \in -K \vspace{2mm}\\
\stackrel{\textmd{by }(\ref{a89n})}{\Longrightarrow} & \frac{f(\bar x+td_\nu)-f(\bar{x})}{t}+C_\nu \frac{x_\nu-\bar{x}}{\|x_\nu-\bar{x}\|}\in -K \vspace{2mm} \\
\hspace*{1mm}\stackrel{\nu\to \infty}{\Longrightarrow} & \frac{f(\bar x+td)-f(\bar{x})}{t} \in -K, \vspace{2mm} \\
\end{array}$$
So, according to the convexity of $\mu^*\circ f$, we have
$$\eta^Td\leq \frac{(\mu^*\circ f)(\bar{x}+td)-(\mu^*\circ f)(\bar{x})}{t} \leq 0,~~\forall\Big(\mu^*\in -K^*,~\eta \in \partial_c (\mu^*\circ f)(\bar x)\Big).$$
Thus, $d\in G_2(\bar x)$, leading to $0 \neq d\in T_\Omega (\bar{x}) \cap G_2(\bar{x}).$ This contradicts the assumption, and the proof of part (ii) is completed.\\
(iii) This part results from parts (i) and (ii) accompanying Definition \ref{cq1}.
\qed


As mentioned before, if $K=\mathbb{R}^p_+,$ the feasible set $\Omega$ is convex, and $f_i,$ $i=1,2,\ldots,p$, are convex, then $G_1(\bar x)=G_2(\bar x),$ i.e. CQ1 automatically holds. This fact leads to the following corollary, derived from Theorem \ref{3part}.

\begin{corollary}\label{cor23}
Assume that $K=\mathbb{R}^p_+,$ the feasible set $\Omega$ is convex, and $f_i,$\linebreak $i=1,2,\ldots,p$, are convex. Then, $\bar{x} \in \Omega$ is a norm-based robust efficient solution of (\ref{prb}) w.r.t. $K=\mathbb{R}^p_+$ if and only if $T_{\Omega}(\bar{x})\cap G_1(\bar{x})=\{0\}$.
\end{corollary}


The following example shows that Theorem \ref{3part}(iii) is not true without~CQ1.

\begin{example}
Consider a VOP $$\min(f_1(x), f_2(x))~~s.t.~~x\in \mathbb{R},$$
with ordering cone $$K=\{ (\mu_1 ,\mu_2)\in \mathbb{R}^2:~ \mu_1 +\mu_2 \geq 0,~ \mu_1 \geq 0\}$$
and $$f_1(x):=\max\{0,x\},~~f_2(x):=\min\{0,-x\}.$$
It is not difficult to see that $f=(f_1, f_2)$ is $K$-convex. Let $\bar{x}=0$. We have
$$\partial_c f(\bar{x})=[0,1]\times [-1,0]$$ and $$-K^*=\{ (\mu_1 ,\mu_2)\in \mathbb{R}^2:~\mu_1\geq \mu_2\geq0\}.$$ So, $G_1(\bar{x})=\{0\}$
and $\partial_c(\mu^*\circ f)(\bar x)=[0,\mu_1-\mu_2]$ for each $(\mu_1,\mu_2)\in -K^*$. Hence,
$G_2(\bar x)=(-\infty,0].$ Therefore, CQ1 is not fulfilled. It is seen that\linebreak $T_\Omega(\bar{x}) \cap G_1(\bar{x}) =\{0\}$ while $\bar{x}$ is not a norm-based robust efficient solution of the considered problem.\qed
\end{example}

The above example has another message that unlike $K=\mathbb{R}_+^p$, for $K\neq \mathbb{R}^p_+$ one can not derive CQ1 from $K$-convexity of the objective function.

The next result characterizes robustness by means
of the normal cone and Clarke's generalized gradient.

\begin{theorem}\label{nec2} Assume that $\Omega$ is convex and $\bar{x} \in \Omega$.
\begin{itemize}
\item[(i)] If $\bar{x}$ is a norm-based robust efficient solution of
(\ref{prb}) w.r.t. $K$, then
$$co\,(-\partial_c f(\bar{x})K^*)+N_\Omega (\bar{x})=\mathbb{R}^n.$$
\item[(ii)] If $f$
is $K$-convex and
\begin{equation}\label{suf11}
pos \bigg(\displaystyle\bigcup_{\mu^* \in -K^*}\partial_c (\mu^*\circ f) (\bar{x})\bigg)+N_\Omega (\bar{x})=\mathbb{R}^n,
\end{equation}
then $\bar{x}$ is a norm-based robust efficient solution of
(\ref{prb}) w.r.t. $K$.
\end{itemize}
\end{theorem}
{\it \textbf{Proof}}
(i) According to Theorem \ref{3part}(i),  norm-based robust efficiency implies
$$T_\Omega (\bar{x}) \cap G_1(\bar{x}) =\{0\}.$$
Now, by \cite[Corollary 16.4.2]{roc}, we get
$$cl\big(N_\Omega (\bar{x}) + G_1(\bar{x})^*\big) =\mathbb{R}^n.$$
On the other hand,
$$G_1(\bar{x})^*= cl\,co(-\partial_c f(\bar{x})K^*).$$ So,
$$\mathbb{R}^n=cl\big(cl\, co(-\partial_c f(\bar{x})K^*)+N_\Omega (\bar{x})\big)=cl\big(co(-\partial_c f(\bar{x})K^*)+N_\Omega (\bar{x})\big).$$
Hence, the closure of the convex set
$co(-\partial_c f(\bar{x})K^*)+N_\Omega (\bar{x})$ coincides with $\mathbb{R}^n$. This leads to
$$co(-\partial_c f(\bar{x})K^*)+N_\Omega (\bar{x})=\mathbb{R}^n.$$
(ii) By setting $$\Gamma_{\bar x} :=pos\big(\displaystyle\bigcup_{\mu^* \in -K^*}\partial_c (\mu^*\circ f) (\bar{x})\big),$$ and applying \cite[Corollary 16.4.2]{roc} on (\ref{suf11}), we have
$$\Gamma_{\bar{x}}^* \bigcap T_\Omega (\bar{x})= \{0\}.$$
On the other hand, $G_2(\bar{x})=\Gamma^*_{\bar{x}}$.
Therefore, $$T_\Omega(\bar{x}) \cap G_2(\bar{x}) =\{0\},$$ and the proof is completed due to Theorem \ref{3part}(ii).
\qed

Corollaries \ref{cobe1}, \ref{cobe12}, and \ref{cobe13} are direct results of Theorems \ref{3part} and \ref{nec2}.

\begin{corollary}\label{cobe1}
Let $\Omega$ be convex, $f_i,~i=1,2,\ldots ,p,$
be convex, and $\bar x\in \Omega$. Then the following three assertions are equivalent.
\begin{itemize}
\item [(i)] $\bar x$ is a norm-based robust efficient solution of (\ref{prb}) w.r.t. $K=\mathbb{R}^p_+$;
\item[(ii)] There exists no $d \in T_\Omega (\bar{x}) \setminus \{0\}$ such that for any $i=1,2,\ldots ,p$
and any $\xi \in \partial f_i(\bar{x})$, $\xi^Td \leq 0$;
\item[(iii)]  $pos\bigg(\displaystyle\bigcup_{i=1}^p\partial f_i (\bar{x})\bigg)+N_\Omega (\bar{x})=\mathbb{R}^n.$
\end{itemize}
\end{corollary}

\begin{corollary}\label{cobe12}
Let $\Omega$ be convex and $f$
be $K$-convex and continuously differentiable at $\bar x\in \Omega$. Then the following three statements are equivalent.
\begin{itemize}
\item [(i)] $\bar x$ is a norm-based robust efficient solution of (\ref{prb}) w.r.t. $K$;
\item[(ii)] There exists no $d \in T_\Omega (\bar{x}) \setminus \{0\}$ such that $\nabla f(\bar{x})^Td\in-K$;
\item[(iii)]  $-\nabla f(\bar{x})^TK^*+N_\Omega (\bar{x})=\mathbb{R}^n.$
\end{itemize}
\end{corollary}

\begin{corollary}\label{cobe13}
Let $\Omega$ be convex and $f_i,~i=1,2,\ldots ,p,$ be convex and continuously differentiable at $\bar{x}\in \Omega$.
Then the following three assertions are equivalent.
\begin{itemize}
\item [(i)] $\bar x$ is a norm-based robust efficient solution of (\ref{prb}) w.r.t. $K=\mathbb{R}^p_+$;
\item[(ii)] There exists no $d \in T_\Omega (\bar{x}) \setminus \{0\}$ satisfying $\nabla f_i(\bar{x})^Td \leq 0$, $i=~1,2,\ldots ,p$;
\item[(iii)]  $pos\big\{\nabla f_1(\bar{x}), \nabla f_2(\bar{x}),\ldots, \nabla f_p(\bar{x})\big\}+N_\Omega (\bar{x})=\mathbb{R}^n.$
\end{itemize}
\end{corollary}



\section{Problems with Conic Constraints}\label{conic}

Consider a VOP with conic constraints as follows:$\vspace{-5mm}$\\
\begin{equation}\label{p28}
\min\, f(x)~~s.t.~~g(x) \in -Q.
\end{equation}
Here, $f: \mathbb{R}^n \rightarrow \mathbb{R}^p$ and $g: \mathbb{R}^n \rightarrow \mathbb{R}^q$ are respectively the objective and
 constraint functions, whose components are assumed to be locally Lipschitz. Furthermore, $Q$ is an ordering cone in $\mathbb{R}^q$ with nonempty interior. The feasible set of (\ref{p28}) is
$$\Omega_1= \{x \in \mathbb{R}^n: g(x)\in -Q\}=\{x \in \mathbb{R}^n:~ (\lambda^* \circ g)(x) \leq 0, ~~\forall \lambda^* \in -Q^*\cap \mathbb{S}\}.$$

Consider a function $G:\mathbb{R}^n \longrightarrow \mathbb{R}$ defined by
$$
G(x):=\displaystyle\max_{\lambda^*}\Big\{(\lambda^* \circ g)(x):~ \lambda^* \in -Q^*\cap \mathbb{S}\Big\},~~~x\in\mathbb{R}^n,
$$
and set
$$I(x):=\Big\{\lambda^* \in -Q^*\cap \mathbb{S}:~G(x)= (\lambda^* \circ g)(x)\Big\},~~~x\in\mathbb{R}^n,$$
$$A(x)=\{\lambda^* \in -Q^*\cap \mathbb{S}:~(\lambda^* \circ g)(x)=0\},~~~x\in\mathbb{R}^n.$$
It is evident that $$\Omega_1=\{x\in\mathbb{R}^n:G(x)\leq 0\}.$$

The following lemma constructs the main result of the current section. Hereafter,
$$\begin{array}{c}
\Upsilon_{\bar{x}} :=\displaystyle\bigcup_{\lambda^* \in A(\bar x)}\partial_c (\lambda^*\circ g) (\bar{x}),\\
\mathcal{D}_{\bar{x}}:=\Big\{ d\in \mathbb{R}^n:~ (\lambda^* \circ g)^{\circ} (\bar{x};d) \leq 0,~~ \forall \lambda^* \in A(\bar{x}) \Big\}.
\end{array}$$
\begin{lemma}\label{lem0928}
If $\bar x\in \Omega_1$ and $0\notin \partial_c G(\bar x)$, then
\begin{itemize}
\item[(i)] $\Big\{ d\in \mathbb{R}^n:~ (\lambda^* \circ g)^{\circ} (\bar{x};d) \leq 0,~ \forall \lambda ^* \in A(\bar{x}) \Big\} \subseteq T_{\Omega_1} (\bar{x});$\vspace{2mm}
\item[(ii)] $N_{\Omega_1} (\bar{x}) \subseteq cl\,pos \bigg(\displaystyle \bigcup_{\lambda^* \in A(\bar x)}\partial_c (\lambda^* \circ g)(\bar x)\bigg);$
\item[(iii)] If $g$ is $Q$-convex, then the inequalities given in (i) and (ii) hold as equality.
\end{itemize}
\end{lemma}
{\it \textbf{Proof}}
(i) If $A(\bar x)=\emptyset$, then there is nothing to prove. So, assume $A(\bar x)\neq \emptyset$. By \cite[Theorems 10.34 and 10.42]{cla}, since $G$ is locally Lipschitz around $\bar x$, we get
\begin{equation}\label{tamam}
\Big\{ d\in \mathbb{R}^n:~ G^{\circ} (\bar{x};d) \leq 0\Big\} \subseteq T_{\Omega_1}^c(\bar{x})\subseteq T_{\Omega_1}(\bar{x}),
\end{equation}
where $T_{\Omega_1}^c(\bar{x})$ stands for the Clarke tangent cone to $\Omega_1$ at $\bar x$. Moreover,\linebreak $G(\bar x)=0$ because $A(\bar x)\neq\emptyset$.

Now, consider $d\in\mathbb{R}^n$ satisfying $(\lambda^*\circ g)^{\circ} (\bar{x};d) \leq 0$ for any $\lambda ^* \in A(\bar{x})$. According to the definition of generalized Clarke's directional derivatives, there exist sequences
$x_\nu \rightarrow \bar{x}$, $t_\nu\downarrow 0$, $\lambda^*_\nu \in I(x_\nu+t_\nu d)$, $\bar\lambda^*_\nu \in I(x_\nu)$ such that
\begin{equation}\label{eq897}
\begin{array}{ll}
 G^{\circ} (\bar{x};d) &=\displaystyle \lim_{\nu\rightarrow \infty} \frac{(\lambda^*_\nu \circ g)(x_\nu+t_\nu d) -(\bar\lambda^*_\nu \circ g)(x_\nu)}{t_\nu}\vspace*{1mm}\\
 &\leq \displaystyle \lim_{\nu\rightarrow \infty} \frac{(\lambda^*_\nu \circ g)(x_\nu+t_\nu d) -(\lambda^*_\nu \circ g)(x_\nu)}{t_\nu}.
\end{array}
\end{equation}
The sequence $\{\lambda_\nu^*\}$ is bounded and, by working with subsequences if necessary, one may assume that this sequence converges to some $\lambda^*\in -Q^*\cap \mathbb{S}$. Furthermore,  as $G$ and $\lambda_\nu^*\circ g$ are continuous at $\bar x$,
$$\begin{array}{ll}
\lambda^*_\nu \in I(x_\nu+t_\nu d),~~\forall\nu&\Longrightarrow G(x_\nu+t_\nu d)=(\lambda_\nu^*\circ g)(x_\nu+t_\nu d),~~\forall\nu\vspace{1mm}\\
&\Longrightarrow 0=G(\bar x)=(\lambda^*\circ g)(\bar x)\Longrightarrow \lambda^*\in A(\bar x).
\end{array}$$
Moreover, due to the locally Lipschitzness of $g_i$ functions, the sequence\linebreak $\left\{\frac{g(x_\nu+t_\nu d) -g(x_\nu)}{t_\nu}\right\}$ is bounded, and hence,
\begin{equation}\label{eq789}
\displaystyle\lim_{\nu\rightarrow \infty} \frac{\left((\lambda^*_\nu-\lambda^*) \circ g\right)(x_\nu+t_\nu d) -(\left(\lambda^*_\nu-\lambda^*) \circ g\right)(x_\nu)}{t_\nu}=0.
\end{equation}
By (\ref{eq897}) and (\ref{eq789}), we get
$$
G^{\circ} (\bar{x};d) \leq \displaystyle \lim_{\nu\rightarrow \infty} \frac{(\lambda^*\circ g)(x_\nu+t_\nu d) -(\lambda^*\circ g)(x_\nu)}{t_\nu}\leq (\lambda^*\circ g)^\circ (\bar x;d)\leq 0.
$$
Hence, $G^{\circ} (\bar{x};d) \leq 0$, and this completes the proof of part (i) due to (\ref{tamam}).

(ii) According to part (i), we have
$N_{\Omega_1} (\bar{x}) =T_{\Omega_1} (\bar{x})^*\subseteq \mathcal{D}_{\bar{x}}^*.$ Furthermore, $\Upsilon_{\bar{x}}^* =\mathcal{D}_{\bar{x}}$, leading to $\mathcal{D}_{\bar{x}}^*=\Upsilon_{\bar{x}}^{**}=cl\, pos (\Upsilon_{\bar{x}})$. So, $N_{\Omega_1} (\bar{x})\hspace*{-.25mm} \subseteq \hspace*{-.25mm} cl\, pos (\Upsilon_{\bar{x}})$.

(iii) As $g$ is $Q$-convex, $G$ is convex. So, due to $0 \not \in \partial_c G(\bar{x})$, there exists $x\in \Omega_1$ such that $(\lambda^* \circ g)(x) <0$ for any $\lambda^* \in -Q^*\cap \mathbb{S}$. On the other hand, the function
$(\lambda^* ,x) \longmapsto \lambda^* \circ g(x)$ is continuous on $(-Q^*\cap \mathbb{S} )\times \mathbb{R}^n$ and
$-Q^*\cap \mathbb{S}$ is compact. Therefore, the inequalities given in parts (i) and (ii) hold as
equality because of \cite[Theorem 7.9]{gob-2} and \cite[Proposition 2.9]{cla}.
\qed


Theorem \ref{nec3} is the main achievement of the current section.

\begin{theorem}\label{nec3}
Let $\bar x\in \Omega_1$.
\begin{itemize}
\item[(i)] If $\Omega_1$ is convex and $\bar x$ is a norm-based robust efficient solution of (\ref{p28}) w.r.t. $K$ with $0\notin \partial_c G(\bar x)$, then
$$co(-\partial_c f(\bar{x})K^*)+pos\bigg(\displaystyle\bigcup_{\lambda^* \in A(\bar x)}\partial_c
(\lambda^* \circ g)(\bar x)\bigg)=\mathbb{R}^n.$$
\item[(ii)] Let $f$ be $K$-convex, $g$ be $Q$-convex and $\bar x\in \Omega_1$. If
$$pos \bigg(\displaystyle\bigcup_{\mu^* \in -K^*}\partial_c (\mu^*\circ f)(\bar x)\bigg)+pos\bigg(\displaystyle\bigcup_{\lambda^* \in A(\bar x)}\partial_c
(\lambda^* \circ g)(\bar x)\bigg)=\mathbb{R}^n,$$
then $\bar{x}$ is a norm-based robust efficient solution of (\ref{p28}) w.r.t. $K$.
\end{itemize}

\end{theorem}
{\it \textbf{Proof}}
(i) Apply Lemma \ref{lem0928}(ii) and Theorem \ref{nec2}.\\
(ii) Considering arbitrary $d\in T_{\Omega_1}(\bar{x})\cap G_2(\bar{x})$, we prove $d=0$, and then norm-based robustness of $\bar x$ comes from Theorem \ref{3part}.
According to the assumption of the theorem, there exist
finite sets $T\subseteq -K^*$ and
$S\subseteq A(\bar{x})$,
Clarke's gradients $\xi_{\mu^*} \in \partial(\mu^*\circ f)(\bar x)~(\mu^* \in T)$ and $\zeta_{\lambda^*} \in
\partial(\lambda^* \circ g)(\bar x)~(\lambda^*\in S)$, and scalars $t_{\mu^*}\geq 0~(\mu^*\in T)$ and
$s_{\lambda^*}\geq 0~(\lambda^*\in S)$
such that
$$d=\sum_{\mu^*\in T} t_{\mu^*}\xi_{\mu^*} +\sum_{\lambda^*\in S} s_{\lambda^*}\zeta_{\lambda^*}.$$
On the other hand, as $d\in G_2(\bar{x})$ and $\xi_{\mu^*} \in \partial (\mu^*\circ f)(\bar x)~(\mu^* \in T)$, we have
$d^T\xi_{\mu^*}\leq 0$.
Also, since $d\in T_{\Omega_1}(\bar{x})$, there exist two sequences
$\{d_\nu\}\subseteq \mathbb{R}^n$ and $\delta_\nu \downarrow 0$ such that $d_\nu \longrightarrow d$ and
$\bar{x} +\delta_\nu d_\nu \in \Omega$; see \cite{cla}.
Therefore, $\zeta_{\lambda^*} \in \partial(\lambda^* \circ g)(\bar x)~(\lambda^*\in S)$ implies
$$d^T_\nu\zeta_{\lambda^*}\leq \frac{(\lambda^*\circ g)(\bar{x}+\delta_\nu d_\nu)-(\lambda^*\circ g)(\bar{x})}{\delta_\nu}\leq 0,~~\forall \nu  \in\mathbb{N}$$
leading to $d^T\zeta_{\lambda^*}\leq 0.$ So, we get $d^Td \leq 0$ which means $d=0$.\qed



The following corollaries are direct consequences of the above theorem (when $K=\mathbb{R}^p_+$ and $Q=\mathbb{R}^q_+$ or
$f_i,~i=1,2,\ldots,p,$ and $g_j,~j=1,2,\ldots,q,$ are continuously differentiable). In these corollaries,
$\Omega_1=\{x\in \mathbb{R}^n :g_j(x)\leq 0,~i=1,2,\ldots, q\}$ and $A(x)=\{j\in \{1,2,\ldots, q\} :g_j(x)=0\}$.
\begin{corollary}\label{cor21} Let $\Omega_1$ be convex and $\bar x\in \Omega_1$.
\begin{itemize}
\item[(i)] If $\bar x$ is a norm-based robust efficient solution of (\ref{p28}) w.r.t. $K=\mathbb{R}^p_+$, and $0\notin  co\Big(\displaystyle\bigcup_{j\in A(\bar{x})} \partial_c g_j(\bar x)\Big)$, then
$$pos\bigg(\displaystyle\bigcup_{i=1}^p\partial_c f_j(\bar x)\bigg)+pos\bigg(\displaystyle\bigcup_{j \in A(\bar x)}\partial_c g_j(\bar x)\bigg)=\mathbb{R}^n.$$
\item[(ii)] Let $f_i,~i=1,2,\ldots, p$ and $g_j,~j=1,2,\ldots, p,$ be convex. If
$$pos\bigg(\displaystyle\bigcup_{i=1}^p\partial_c f_j(\bar x)\bigg)+pos\bigg(\displaystyle\bigcup_{j \in A(\bar x)}\partial_c g_j(\bar x)\bigg)=\mathbb{R}^n,$$
then $\bar{x}$ is a norm-based robust efficient solution of (\ref{p28}) w.r.t. $\mathbb{R}^p_+$.
\end{itemize}
\end{corollary}


\begin{corollary}\label{cor22} Let $\Omega_1$ be convex and $f_i,~i=1,2,\ldots ,p,$ and $g_j,~j\in A(\bar{x}),$ be continuously differentiable at $\bar{x}\in \Omega_1$.
\begin{itemize}
\item[(i)] If $\bar x$ is a norm-based robust efficient solution of (\ref{p28}) w.r.t. $K=\mathbb{R}^p_+$, and $0\notin  co\{\nabla g_j(\bar{x}) : j\in A(\bar{x})\}$, then
$$pos\big\{\nabla f_1(\bar{x}), \nabla f_2(\bar{x}),\ldots, \nabla f_p(\bar{x})\big\}+
pos\big\{\nabla g_j(\bar{x}) : j\in A(\bar{x})\big\}=\mathbb{R}^n.$$
\item[(ii)] Let $f_i,~i=1,2,\ldots ,p,$ and $g_j,~j=1,2,\ldots ,q,$ be convex. If
$$pos\big\{\nabla f_1(\bar{x}), \nabla f_2(\bar{x}),\ldots, \nabla f_p(\bar{x})\big\}+
pos\big\{\nabla g_j(\bar{x}) : j\in A(\bar{x})\big\}=\mathbb{R}^n,$$
then $\bar{x}$ is a norm-based robust efficient solution of (\ref{p28}) w.r.t. $\mathbb{R}^p_+$.
\end{itemize}
\end{corollary}

In the following, we investigate the robustness for a semi-infinite VOP.
Consider the following semi-infinite VOP:
\begin{equation}\label{p3}
\min\, f(x)~~s.t.~~g_j(x)\leq 0, ~~~\forall j \in J.
\end{equation}
Here, $f: \mathbb{R}^n \rightarrow \mathbb{R}^p$ is the objective function
(i.e. $f(x) = (f_1(x), \ldots, f_p(x))$, $g_j: \mathbb{R}^n \rightarrow \mathbb{R}~(j\in J)$ are
the constraint functions and $J$ is an infinite index set. We set $\Omega_1$ as the feasible set of (\ref{p3}), i.e.,
$$
\Omega_1 =\{x \in \mathbb{R}^n ~;~ g_j(x) \leq 0, ~~\forall j \in J\}.
$$
Let $f_i~(i=1, \ldots , p)$ and $g_j~(j\in J)$ be locally Lipschitz. Notice that Problem (\ref{p3}) is a special case of (\ref{p28}). We say that the Slater constraint qualification (SCQ) holds for (\ref{p3}) if the following conditions are together satisfied:
\begin{itemize}
\item[(i)] J is compact,
\item[(ii)] The function $(j,x) \longmapsto g_j(x)$ is continuous on $J\times \mathbb{R}^n$,
\item[(iii)] There is a $x^\circ \in \mathbb{R}^n$ such that $g_j(x^\circ) <0,$ for any $j\in J$.
\end{itemize}

Set $J(\bar x)=\{j \in J:~g_j(\bar x)=0\}.$ In the following theorem, we provide a characterization for norm-based robust efficient solutions of $(\ref{p3})$.

\begin{theorem}\label{nec4}
Let $g_j,~j\in J,$ be convex and $\bar{x}\in \Omega_1$.
\begin{itemize}
\item[(i)] Assume that $0\notin co\big\{\bigcup_{j \in A(\bar{x})} \partial g_j(\bar{x})\big\}$ and \textnormal{SCQ} holds for (\ref{p3}). If $\bar x$ is a norm-based robust efficient solution of (\ref{p3}) w.r.t. $K$, then
$$co\big(-\partial_c f(\bar{x})K^*\big)+pos\bigg(\displaystyle\bigcup_{j \in A(\bar x)}\partial
g_j(\bar x)\bigg)=\mathbb{R}^n.$$
\item[(ii)] Assume that $f$ is $K$-convex. Then, $\bar{x}$ is a norm-based robust efficient solution of (\ref{p3}) w.r.t. $K$, if
$$pos \bigg(\displaystyle\bigcup_{\mu^* \in -K^*}\partial_c (\mu^*\circ f)(\bar x)\bigg)+pos\bigg(\displaystyle\bigcup_{j \in A(\bar x)}\partial g_j(\bar x)\bigg)=\mathbb{R}^n.$$
\end{itemize}
\end{theorem}
{\it \textbf{Proof}}
Apply Theorem \ref{nec2} and \cite[Lemma 7.7 and Theorem 7.9]{gob-2}.
\qed

\section{Robustness and Gap Function}\label{gap}

Gap function is one of the important tools for characterizing optimality/efficiency in optimization.
To the best of our knowledge, gap function for multiobjective problems was first developed by Chen et al. \cite{Chen1998}. 
Here, we define a gap function for VOP (\ref{prb}) as a set-valued mapping $\Phi_{gap}:\Omega \times \partial_c f(\cdot ) \rightrightarrows \mathbb{R}^p$ defined as
$$\Phi_{gap}(x,\xi):=\{\xi^T (x-\bar{y}) : \bar{y}\in E_{x,\xi}\},$$
for $x\in \Omega$ and $\xi \in \partial_c f(x)$. In this gap function, $E_{x,\xi}$ is the set of efficient solutions of the following VOP w.r.t. $K$,
$$\max_{y\in \Omega } \, \, \xi^T (x-y).$$


The vector-valued function $f:\mathbb{R}^n \longrightarrow \mathbb{R}^p$ is called $K$-regular at $\bar{x}\in \mathbb{R}^n$, if for any $\mu^* \in -K^*$ the function $\mu^* \circ f$ is regular at $\bar{x}$. If $f$ is $K$-convex, then it is $K$-regular. Also, if $K=\mathbb{R}^p_+$, then $f$ is $K$-regular if and only if $f_i$ is regular for $i=1, 2,\ldots, p$.

The following theorem provides a necessary condition for norm-based robust efficiency.
\begin{theorem}
Let $\Omega$ be convex and $f$ be $K$-regular at $\bar{x}\in \Omega$. If $\bar{x}$ is a norm-based robust efficient solution of (\ref{prb}) w.r.t. $K$ and (\ref{condition}) is fulfilled, then there exists
$\bar{\xi} \in \partial_c f(\bar{x})$ such that $0\in \Phi_{gap}(\bar{x},\bar{\xi})$.
\end{theorem}
{\it \textbf{Proof}}
As $\bar{x}$ is a norm-based robust efficient solution, according to Theorem \ref{nec2}, we get
$$co\big(-\partial_c f(\bar{x})K^{*\circ }\big)+N_\Omega (\bar{x})=\mathbb{R}^n.$$
So, there exist $m\in \mathbb{N}$, $\xi_\nu \in \partial_c f(\bar{x})$, $\mu_\nu^* \in -K^{*\circ }$, and
$t_\nu \geq 0,~\nu =1,2,\ldots ,m,$ such that
$-\sum_{\nu =1}^m t_\nu \xi_\nu \mu_\nu^* \in N_\Omega (\bar{x})$ and $\sum_{\nu =1}^m t_\nu=1$. Now, by (\ref{condition}) and \linebreak$K$-regularity of $f$, we get
$$\begin{array}{ll}
&\xi_\nu \mu_\nu^* \in \partial_c f(\bar{x})\mu^*_\nu =\partial_c (\mu^*_\nu \circ f) (\bar{x}),
~~\forall ~\nu =1,2,\ldots ,m\vspace*{1mm}\\
\Longrightarrow & \displaystyle \sum_{\nu =1}^m t_\nu \xi_\nu \mu_\nu^* \in \sum_{\nu =1}^m t_\nu \partial_c (\mu_\nu ^* \circ f) (\bar{x})= \partial_c ((\sum_{\nu =1}^m t_\nu \mu_\nu^*) \circ f) (\bar{x}) = \partial_c f (\bar{x}) \sum_{\nu =1}^m t_\nu \mu_\nu^{*}.
\end{array}$$
Therefore, there exists $\xi \in \partial_c f(\bar{x})$ such that
$$-\sum_{\nu =1}^m t_\nu \xi \mu_\nu^* =-\sum_{\nu =1}^m t_\nu \xi_\nu \mu_\nu^* \in N_\Omega (\bar{x}).$$
As $\Omega$ is convex,
\begin{equation}\label{fg}
-\displaystyle\sum_{\nu =1}^m t_\nu \mu_\nu^{*T}\xi^T (y -\bar{x})\leq 0,~~ \forall y\in \Omega.
\end{equation}
On the other hand, if $0\not \in \Phi_{gap}(\bar{x},\xi)$, then $\bar{x} \not \in E_{\bar{x},\xi}$ and so, there exists $y \in \Omega$ such that $\xi^T(\bar{x}-y)\in K\setminus \{0\}$. This leads to
$\mu_\nu^{*T} \xi^T (y-\bar{x})< 0$ for each $\nu =1,2,\ldots ,m$, due to $\mu_\nu^* \in -K^{*\circ}$. Hence,
$-\displaystyle \sum_{\nu =1}^m t_\nu \mu_\nu^{*T} \xi^T (y-\bar{x})> 0.$
This strict inequality contradicts (\ref{fg}) and the proof is complete.
\qed


\begin{corollary}
Let $\Omega$ be convex and $f_i,~i=1,2,\ldots ,p,$ be continuously differentiable. If $\bar{x} \in \Omega$ is a norm-based robust efficient solution of (\ref{prb}) w.r.t. $K$, then $0\in \Phi_{gap}(\bar{x},\nabla f(\bar{x}))$.
\end{corollary}

\begin{corollary}
Let $\Omega$ be convex and $f_i,~i=1,2,\ldots,p,$ be regular. If $\bar{x} \in \Omega$ is a norm-based robust efficient solution of (\ref{prb}) w.r.t. $K=\mathbb{R}^p_+$, then there exists
$\bar{\xi} \in \partial_c f(\bar{x})$ such that $0\in \Phi_{gap}(\bar{x},\bar{\xi})$.
\end{corollary}



\section{Conclusions}\label{Conclusion}

Investigation and characterization of the norm-based robust solutions of VOPs was the main aim of the current work. After addressing some basic notions, we obtained necessary and sufficient conditions for norm-based robustness utilizing
two new directions, defined invoking Clarke's generalized gradient. Furthermore, we developed these conditions for
norm-based robust efficient solutions of VOPs with conic constraints and semi-infinite VOPs.
Moreover, we derived a necessary condition for norm-based robustness by means of a nonsmooth gap function. In addition to general results, we analysed the problem for special cases.



%


\begin{thebibliography}{}


\bibitem{ben-1}
Ben-Tal, A., Ghaoui, L. El., Nemirovski, A.: Robust Optimization. Princeton Ser. Appl.
Math., Princeton University Press, NJ (2009)

\bibitem{ben-2}
Ben-Tal, A., Nemirovski, A.: Robust optimization methodology and applications,
Mathematical Programming. 92, 453-480 (2002)

\bibitem{ber}
Bertsimas, D., Brown, D. B., Caramanis, C.: Theory and applications of robust optimization,
SIAM Review. 53, 464-501 (2011)

\bibitem{ehr-2}
Ehrgott, M., Ide, J., Sch\"{o}bel, A.: A Minmax robustness for multi-objective optimization problems,
European Journal of Operational Research. 239, 17-31 (2014)

\bibitem{fli-2}
Fliege, J.,  Werner, R.: Robust multiobjective optimization  applications in portfolio optimization,
European Journal of Operational Research. 234, 422-433 (2014)

\bibitem{geo}
Georgiev, P. Gr., Luc, D. T., Pardalos, P.: Robust aspects of solutions in deterministic multiple objective linear programming, European Journal of Operational Research. 229, 29-36 (2013)

\bibitem{gob-1}
Goberna, M. A., Jeyakumar, V. Li. G., L\'{o}pez, M. A.: Robust solutions to multi-objective linear programs with uncertain data, European Journal of Operational Research. 242, 730-743 (2015)

\bibitem{ide-1}
Ide, J., K$\ddot{o}$bis, E.: Concepts of efficiency for uncertain multi-objective optimization problems based on set order relations, Mathematical Methods of Operations Research. 80, 99-127 (2014)

\bibitem{ide-3}
Ide, J., Sch$\ddot{o}$bel, A.: Robustness for uncertain multi-objective optimization: a survey and analysis of different concepts, OR Spectrum. 38, 235-271 (2016)

\bibitem{kur}
Kuroiwa, D.,  Lee, G. M.: On robust multi-objective optimization, Vietnam Journal of Mathematics. 40, 305-317 (2012)

\bibitem{zam}
Zamani, M., Soleimani-damaneh, M., Kabgani, A.: Robustness in nonsmooth nonlinear multi-objective programming, European Journal of Operational Research. 247, 370-378 (2015)

\bibitem{our-new}
Rahimi, M., Soleimani-damaneh, M.: Robustness in deterministic vector optimization. Journal of Optimization Theory and Applications 179, 137-162 (2018)

\bibitem{roc}
Rockafellar, R. T.: Convex Analysis. Second Printing, Princeton University, Princeton, New Jersey (1972)

\bibitem{cla}
Clarke, F.: Functional Analysis, Calculus of Variations and Optimal Control. Springer-Verlag, London (2013)

\bibitem{gob-2}
Goberna, M. A., L\'{o}pez, M. A.: Linear Semi-infinite Optimization. John Wiley \& Sons, Berlin (1998)

\bibitem{Chen1998}
Chen, G.Y., Goh,  V.J., Yang, X.Q.: The gap function of a convex multicriteria optimization problem. European Journal of Operational Research, 111, 142-151 (1998)

\end{thebibliography}

\end{document}